\documentclass[11pt]{amsart}
\usepackage{geometry}
\usepackage[all]{xy}
\usepackage{graphicx}
\usepackage{amssymb}
\usepackage{amsfonts}
\usepackage{epstopdf}
\usepackage{amsmath, amsthm}
\usepackage{mathrsfs}
\usepackage{amscd,amsthm,amsmath,amsgen,amsfonts, amssymb, latexsym, epsfig, euler, epic, graphicx, psfrag}

\newcommand{\bb}[1]{\mathbb{#1}}
\newcommand{\lie}[1]{\mathfrak{#1}}

\def\inv{^{-1}}

\theoremstyle{plain}
\newtheorem{theorem}{Theorem}[section]

\theoremstyle{definition}
\newtheorem{example}[theorem]{Example}
\newtheorem{remark}[theorem]{Remark}

\begin{document}
\title{Localization and Specialization for Hamiltonian Torus Actions.}

\author{Milena Pabiniak}

\address{Milena Pabiniak,
CAMGSD, Departamento de Matem\'atica,
Instituto Superior T\'ecnico, Lisboa, Portugal}
\email{mpabiniak@math.ist.utl.pt}

\thanks{\today}

\begin{abstract}
We consider a Hamiltonian action of $n$-dimensional torus, $T^n$, on a compact symplectic manifold $(M,\omega)$ with $d$ isolated fixed points. 
For every fixed point $p$ there exists (though not unique) a class $a_p \in H^*_{T}(M; \bb{Q})$ such that the collection $\{a_p\}$, over all fixed points, 
forms a basis for $H^*_{T}(M; \bb{Q})$ as an $H^*(BT; \bb{Q})$ module. 
The map induced by the inclusion, $\iota^*:H^*_{T}(M; \bb{Q})  \rightarrow H^*_{T}(M^{T}; \bb{Q})= \oplus_{j=1}^{d}\bb{Q}[x_1, \ldots, x_n] $ is injective. 
We use such classes $\{a_p\}$ to give necessary and sufficient conditions for $f=(f_1, \ldots ,f_d)$ in $\oplus_{j=1}^{d}\bb{Q}[x_1, \ldots, x_n]$ 
to be in the image of $\iota^*$, i.e. to represent an equiviariant cohomology class on $M$. In the case when $T$ is a circle and present these conditions explicitly. We explain how to combine this $1$-dimensional solution with Chang-Skjelbred Lemma in order to obtain the result for a torus $T$ of any dimension. Moreover, for a GKM $T$-manifold $M$ our techniques give combinatorial description of $H^*_{K}(M; \bb{Q})$, for a generic subgroup $K \hookrightarrow T$, 
even if $M$ is not a GKM $K$-manifold.
\end{abstract}

\maketitle

\tableofcontents

\section{Introduction}

Suppose that a compact Lie group $G$ acts on a compact, closed, connected and oriented manifold $M$. 
The equivariant cohomology ring $H^*_{G}(M;R):=H^*(M \times_G EG;R)$, with coefficients in a ring $R$, 
encodes topological information about the manifold and the action. In the case of a Hamiltonian action on a 
symplectic manifold, a variety of techniques has made computing $H^*_{G}(M;R)$ tractable. 
The work of Goresky-Kottwitz-MacPherson \cite{GKM} describes this ring combinatorially when $G$ is a torus, $R$  a field, 
and the action has a very specific form. We give a more general description that has a similar flavor. 
A theorem of Kirwan \cite{KirwanCohomology} states that the inclusion of the fixed points induces an injective map in equivariant cohomology. We quote this result below, following Tolman and Weitsman \cite{TW2}.
\begin{theorem}[Kirwan, \cite{KirwanCohomology}]\label{KirwanInjectivity} Let a torus $T$ act on a symplectic compact connected manifold $(M,\omega)$ 
in a Hamiltonian fashion and let $\iota: M^{T}
 \rightarrow M$ denote the natural inclusion of fixed points. Then the induced map $\iota^*:H^*_{T}(M; \bb{Q})  \rightarrow H^*_{T}(M^{T}; \bb{Q}) $ 
is injective. 
 If $M^T$ consists of isolated points then also $\iota^*:H^*_{T}(M; \bb{Z})  \rightarrow H^*_{T}(M^{T};
  \bb{Z}) $ is injective.
\end{theorem} 
 \noindent
If there are $d$ fixed points then $H^*_{T}(M^{T}; \bb{Q})= \oplus_{j=1}^{d} \bb{Q}[x_1, \ldots, x_n]$, where $n$ is the dimension of the torus.
Therefore we can think about an equivariant cohomology class in $H^*_T(M^T;\bb{Q})$ as a $d$-tuple of polynomials $f=(f_1, \ldots ,f_d)$, 
with each $f_j$ in $\bb{Q}[x_1, \ldots, x_n]$. The goal of this paper is to give necessary and sufficient conditions for a $d$-tuple of polynomials 
to be in the image of $\iota^*$, that is to represent an equiviariant cohomology class on $M$. 
\\Notation: By abuse of language we
will often say that a $d$-tuple of polynomials $f=(f_1, \ldots ,f_d)$ ``is'' or ``represents'' an equivariant cohomology class if it
is the image under $\iota^*$ of an honest (unique) equivariant cohomology class
on $M$.
\begin{remark}\label{reductiontoonedim} {\bf Reducing the problem to the case $T=S^1$.} The following result of Chang and Skjelbred \cite{ChangS} guarantees that we only need to consider the case of an $S^1$ action. Example \ref{CSusage} shows how to combine Theorems \ref{main} and \ref{CSLemma} in order to obtain information about $H^*_T(M; \bb{Q})$ for a torus $T$ of higher dimension. 
\begin{theorem}[Chang, Skjelbred, \cite{ChangS}]\label{CSLemma}
The image of $\iota^*:H^*_{T}(M; \bb{Q})  \rightarrow H^*_T(M^T; \bb{Q}) $ is the set
$$ \bigcap_H \iota^*_{M^H}(H^*_{T}(M^H; \bb{Q})),$$
where the intersection in $H^*_T(M^T; \bb{Q})$ is taken over all codimension-one subtori $H$ of $T$, and $\iota_{M^H}$ is the inclusion of $M^T$ into $M^H$.
\end{theorem}

In fact the only nontrivial contributions to this intersection are those codimension $1$ subtori $H$ which appear as isotropy groups of some elements of $M$ 
(that is $M^H \neq M^{T}$).
\end{remark}

Therefore we will consider
 a circle $S^1$ acting on a 
compact, connected and closed 
symplectic manifold $(M^{2n},\omega)$ in a 
Hamiltonian fashion with isolated fixed points and moment map $\mu:M \rightarrow \bb{R}$. 
Unless otherwise stated, all the manifold considered in this paper are assumed to be compact, closed and connected.
It turns out that with these assumptions we are in the Morse Theory setting. Averaging the symplectic form if necessary we can assume that $\omega$ is $S^1$-invariant. Such $(M, \omega)$ can be equipped with an invariant compatible almost complex structure $J$ (compatible means that $\omega( ., J .)$ is a Riemannian metric). Then for any fixed point $p$ there are integers $\{\eta_1, \ldots, \eta_n\}$ such that the $S^1$ action on $T_pM$ is isomorphic to the $S^1$ action on $\bb{C}^n$ with weights $\{\eta_1, \ldots, \eta_n\}$. These integers are called the (isotropy) weights of the $S^1$ action at a fixed point $p$. Since the set of compatible almost complex structures is contractible, the set of weights of the $S^1$ representation on $T_pM$ does not depend on $J$.
\begin{theorem}[Frankel \cite{Frankel}, Kirwan \cite{KirwanCohomology}]\label{mmapmorse}
In the above setting, the moment map $\mu$ is a perfect Morse function on $M$ (for both ordinary and equivariant cohomology). 
The critical points of $\mu$ are the fixed points of $M$, and the index of a critical point $p$ is precisely 
twice the number of negative weights of the circle action on $T_pM$.
\end{theorem}
The Morse function is called \textbf{perfect} if the number of critical points of index $k$ is equal to the dimension of $k$-th cohomology group.
The action of a torus of higher dimension also carries a Morse function.
For $\xi \in \lie{t}$ we define $\Phi^{\xi}:M \rightarrow \bb{R}$, the component of moment map along $\xi$, by $\Phi^{\xi}(p)= \langle \Phi,\xi \rangle$.
We call $\xi \in \lie{t}$ {\bf generic}
if $\langle \eta, \xi \rangle \neq 0$ for each weight $\eta \in \lie{t}^*$ of $T$ action on
$T_p M$, for  every $p$ in the fixed points set $M^T$.
For a generic, rational $\xi$, 
$ \Phi^\xi$ is a Morse function with critical set $M^T$. This map is a moment map for the action of a subcircle $S \hookrightarrow T$ generated by $\xi \in \lie{t}$.
Using Morse Theory, Kirwan constructed equivariant cohomology classes that form a basis for integral equivariant cohomology ring of $M$. 
Then the existence of a basis for rational equivariant cohomology ring of $M$ follows. 
We quote this theorem with the integral coeficients, and action of a torus $T$ of any dimension, although in this paper we work mostly with rational coefficients and circle actions.
\begin{theorem}[Kirwan, \cite{KirwanCohomology}]\label{kirwanclasses}
Let a torus $T$ act on a symplectic compact manifold $M$ with isolated fixed points,
and let $\mu = \Phi^\xi:M \rightarrow \bb{R}$ be a component of moment map $\Phi$ along generic $\xi \in \lie{t}$.
Let $p$ be any fixed point of index $2k$ and let $w_1, \ldots, w_k$ be the negative weights of the $T$ action on $T_pM$. 
Then there exists a class $a_p \in H^{2k}_{T}(M; \bb{Z})$ such that
\begin{itemize}
\item $a_p|_p=  \Pi_{i=1}^k\,w_i$; 
\item $a_p|_{p'}=0$ for all fixed points $p' \in M^T \setminus \{p\}$ such that $\mu(p') \leq \mu(p)$.
\end{itemize}   
Moreover, taken together over all fixed points, these classes are a basis for the cohomology $H^*_{T}(M; \bb{Z})$ as an $H^*(BT; \bb{Z})$ module.
\end{theorem}
\noindent
In the above theorem we use the convention that the empty product is equal to $1$.
We will call the above classes \textbf{Kirwan classes}.
These classes may be not unique.
Goldin and Tolman use a different basis for the cohomology ring $H^*_T(M; \bb{Z})$ in [GT]. 
They additionally require $a_p|_{p'}=0$ for all fixed points $p' \neq p$ of index less then or equal $2k$ (where $2k$ is the index of $p$). 
Goldin and Tolman's classes, if they exist, are unique. Therefore they are called \textbf{canonical classes}. 
For our purposes, it is enough to have some basis for the rational equivariant cohomology ring with respect to the given circle action, and with the following property
\begin{itemize}
 \item[($\star$)] elements of the basis are in such a bijection with the fixed points that for a class $\alpha$ corresponding to 
a fixed point of index $2k$, have that if $\iota^*(\alpha)=(f_1,\ldots,f_d) \in  \oplus_{j=1}^d \,\bb{Q}[x]=H^*_{S^1}(M^{S^1}) $ then each $f_j$ is $0$ or a homogeneous polynomial of degree $k$. \end{itemize}
We will call a basis satisfying condition ($\star$) \textbf{a basis of generating classes}. Kirwan classes and Goldin-Tolman canonical classes satisfy the above condition.  

In this paper, we show how to obtain relations describing the image of $\iota^* (H_{S^1}^*(M)) \subset H_{S^1}^*(M^{S^1})$.
For a fixed point $p$ let $e(p)$ be the image of the equivariant Euler class of the tangent bundle evaluated at $p$, which in this case is equal to the product of weights of the circle action.
The Main Theorem is:
\begin{theorem}\label{main}
Let a circle act on a closed compact connected symplectic manifold $M$ in a Hamiltonian fashion, with isolated fixed points $p_1, \ldots, p_d$. Suppose we are given a basis  $\{a_p\}$ of $H_{S^1}^*(M; \bb{Q})$, satisfying condition ($\star$). Let $f=(f_1, \ldots ,f_d) \in \oplus_{j=1}^d \,\bb{Q}[x]=H^*_{S^1}(M^{S^1}; \bb{Q})$. 
Then $f$ is an image (under $\iota^*$ from Theorem \ref{KirwanInjectivity}) of an equivariant cohomology class on $M$ if and only if for every fixed point $p$ 
of index $2k$, $0 \leq k < n$ 
we have
\begin{equation}\label{relations}
 \sum_{j=1}^d \frac{f_j\,a_p(p_j)}{e(p_j)} \in \bb{Q}[x],\end{equation}
where $a_p(p_j)$ denotes $\iota_{p_j}^*(a_p)$, with $\iota_{p_j}:p_j \hookrightarrow M$ the inclusion of fixed point $p_j$ into $M$, that is $\iota^*(a_p)=(a_p(p_1), \ldots, a_p(p_d))\in \oplus \bb{Q}[x]$.
\end{theorem} 
\noindent
Note that if $p$ is a fixed point of index $2n$, this condition is automatically satisfied.
This is because $a_p$ is nonzero only at $p$, and there its value is the Euler class $e(p)$. 
Therefore it is sufficient to check the above condition only for points of index strictly less then $2n=\textrm{dim}M$.

An important ingredient of the proof is the Atiyah-Bott, Berline-Vergne (ABBV) localization theorem. 
  
\begin{theorem}[ABBV Localization, \cite{AB}\cite{BV}]\label{ABBV}
Let $M$ be a compact connected oriented manifold equipped with an $S^1$ action with isolated fixed points, and let $\alpha \in H^*_{S^1}(M; \bb{Q})$. 
Then as elements of $H^*(BS^{1};\bb{Q})=\bb{Q}[x]$,
$$ \int_M \alpha = \sum_p \frac{\alpha|_p}{e(p)},$$
where the sum is taken over all the fixed points.
\end{theorem}
\begin{remark}
If $f$ is a cohomology class, then so is $f\cdot a_p$. Applying the Localization Theorem to the class $f\cdot a_p$ we see that these conditions must be satisfied.
The interesting part of the theorem is that they are sufficient to describe $H_{S^1}^*(M)$ as a subring of $H_{S^1}^*(M^{S^1})$.
\end{remark}
\begin{remark}\label{sphere} \textbf{Connection with the GKM Theorem.}
We now recall the GKM Theorem and therefore for a moment we work with a torus of any dimension.
Let $M$ be a compact, connected, symplectic manifold with a Hamiltonian, effective action of a torus $T=T^n$ and with finitely many fixed points $\{p_1,\ldots, p_d\}$.
Let $N\subset M$ be the set of points whose orbits under the $T$ action are 1-dimensional. The \textbf{one-skeleton} of $M$ is the closure $\overline{N}$.
The manifold $M$ is called 
a \textbf{GKM manifold} if $N$ has finitely many connected components $N_\alpha$. Then for each such component $N_\alpha$ its closure $\overline{N}_\alpha$ is diffeomorphic to a sphere fixed by a codeminesion one subtorus $T_{\alpha}$ of $T$, with residual circle acting by rotation with some weight $w_{\alpha}$, and fixing the north and south poles, $n_{\alpha}, s_{\alpha} \in M^T$. For any class $f \in H^*_T(M^T)$ let $f_{|p_j} $ denote its restriction to fixed point $p_j$.
\begin{theorem}[\cite{GKM},\cite{TW2}]\label{GKMtheorem}
Let $M$ be a GKM manifold with a Hamiltonian torus
action by $T$. Let $M^T$ be the fixed point set, and $\overline{N}$ be the one-skeleton.
Let $\iota:M^T\hookrightarrow M$ be the inclusion of the
fixed point set to $M$ and $j: M^T \hookrightarrow \overline{N}$ be
the inclusion to $\overline{N}$.
The induced maps $\iota^*:H_T^*(M)\rightarrow
H_T^*(M^T)$ and $j^*:H_T^*(\overline{N})\rightarrow H_T^*(M^T)$ on
equivariant cohomology have the same image. 
\end{theorem}
One can extract from the above theorem an explicit describtion of $H^*_T(M)$, namely, $(f_{|p_1}, \ldots, f_{|p_d})\in H^*_T(M^T)$ is in the image of $\iota^*\colon  H^*_T(M) \rightarrow H^*_T(M^T)$ if and only if 
$$\forall_{\overline{N}_{\alpha}}\,\,\,\,f_{|n_{\alpha}}\equiv f_{|s_{\alpha}} \textrm{  mod }w_{\alpha}.$$
The above relations are often called the `` GKM relations''. 

Consider the standard Hamiltonian $S^1$ action on $S^2$ by rotation with a weight $ax$. 
The isolated fixed points are south and north poles which we will denote by $p_1$ and $p_2$ respectively. 
The Goldin-Tolman class associated to $p_1$ is $1$. Theorem \ref{main} says that $f=(f_1,f_2)$ represents an quivariant cohomology class if and only if
$$ \frac{f_1\,a_1(p_1)}{e(p_1)} + \frac{f_2\,a_1(p_2)}{e(p_2)} =\frac{f_1}{ax}+ \frac{f_2}{-ax} = \frac{f_1-f_2}{ax}\in \bb{Q}[x].$$
The above condition is exactly the same as the condition (1) in \cite{GH}. Using the solution for this special case, together with the Chang-Skjelbred Lemma,
Goldin and Holm recover the GKM Theorem in Section 1 and 2 of \cite{GH}.
\end{remark}

Theorem \ref{main} is useful only if we know some basis of generating classes (whose existence is guaranteed by Theorem \ref{kirwanclasses}) and its image under $\iota^*$. Although we cannot compute these classes in general, there are algorithms that work for a wide class of spaces, for example GKM spaces, which
include symplectic toric manifolds and flag manifolds (see \cite{T}). 
For the sake of completeness we will describe an algorithm for obtaining Kirwan classes for symplectic toric manifolds in Appendix \ref{generatorsfortoric}. 
 The choice of $a_p$ assigned to fixed point $p$ may be not unique, even for symplectic toric manifolds.
In the case when moment map is so called ``index increasing''
and the manifold is a GKM manifold, 
uniqueness was proved by Goldin and Tolman in \cite{GT}. 

\begin{remark} {\bf Specialization.} A particularly interesting application of our theorem is when we want to restrict the action of $T$ to an action of a subtorus $S\hookrightarrow T$ such that $M^S=M^T$, 
and compute $\iota^*(H^*_S(M))\subseteq H^*_S(M^S)=H^*_S(M^T)$. 
We call this process \textbf{specialization} of the $T$ action to the action of a subtorus $S$. 
GKM relations are sufficient to describe the image of $  H^*_T(M)$ in $H^*_T(M^T)$, but their ``projections'' are not sufficient to describe the image of $  H^*_S(M)$  
 in $H^*_S(M^T)$.\begin{figure}[h]
	\centering
		\includegraphics[width=0.65\textwidth]{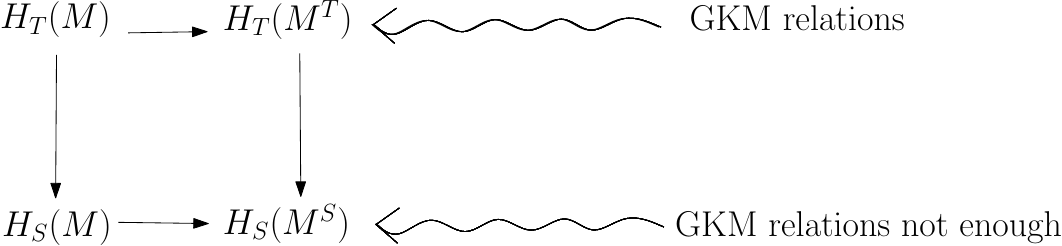}
	\label{fig:GKMRelations2}
\end{figure}
However having generating classes for the $T$ action we can easily compute generating classes for the $S$ action using the projection $\lie{t}^* \rightarrow \lie{s}^*$ (see Appendix \ref{generatorsfortoric} and explicit calculations in Section \ref{examples}, Examples).  
Then the application of Theorem \ref{main} gives precise relations that cut out the image $\iota^*(H^*_S(M))\subseteq H^*_S(M^S)$. 

In particular we can use this method to restrict the torus action on a symplectic toric manifold to a generic circle, i.e. such a circle $S$ for which $M^S=M^T$ (see Examples \ref{toric1} and \ref{toric2}). 
A priori we only require that $M^S$ is finite as we still want to describe $H_S^*(M)$ by analyzing the relations on polynomials defining 
the image $\iota^*(H^*_S(M))\subseteq H^*_S(M^S)=\oplus \bb{Q}[x]$. However it turns out that this requirement implies $M^T=M^S$. 
We can explain this fact using Morse theory. If $\Phi: M \rightarrow \lie{t}^*$ is a moment map for the $T$ action and $\xi \in \lie{t}$ is generic, then $\Phi^{\xi}$, a component of $\Phi$ along $\xi$, is a perfect Morse function with critical set $M^T$. Therefore $\sum \, \textrm{dim}\,H^i(M)=|M^T|$. Similarly, taking $\mu =pr_{\lie{s^*}} \circ \Phi$ for a moment map for the $S$ action, and any generic $\eta \in \lie{s}$, we obtain $\mu^{\eta}$ which 
is also a perfect Morse function for $M$. Thus $|M^S|=\sum \, \textrm{dim}\,H^i(M)=|M^T|$. As obviously $M^T \subset M^S$, the sets must actually be equal.
\end{remark}

The GKM Theorem is a very powerful tool that allows us to compute the image $\iota^* \colon H^*_T(M) \hookrightarrow H^*_T(M^T)$. However this theorem cannot be applied if for some codimension $1$
subtorus $H \hookrightarrow T$ we have dim $M^H >2$. Goldin and Holm in [GH] provide a  generalization of this result to the case where dim $M^H \leq 4$ for all
codimension $1$ subtori $H \hookrightarrow T$. 
An important corollary is that, in the case of Hamiltonian circle actions, with isolated fixed points, on manifolds of dimension $2$ or $4$,
the rational equivariant cohomology ring can be computed solely from the weights of the circle action at the fixed points.
In dimension $2$ this is given for example by the GKM Theorem. In dimension $4$ one can apply the algorithm presented by Goldin and Holm in \cite{GH}
or use the fact that any such $S^1$ action is actually a specialization of a toric $T^2$ action (see \cite{Karshon4mfd}).
If one wishes to compute the integral equivariant cohomology ring, one will need an additional piece of information, so called ``isotropy skeleton`` (\cite{GO}). Godinho in \cite{GO} presents an algorithm for such computation. Information encoded in the isotropy skeleton is essential. There cannot exist an algorithm computing the integral equivariant cohomology only from the fixed point data. Karshon in \cite{KarshonHamiltonian}(Example 1), constructs two $4$-dimensional $S^1$ spaces with the same weights at the fixed point but different integral equivariant cohomology ring. This suggests that we probably should not hope for an algorithm computing the rational equivariant cohomology ring from the weights at the fixed points for manifolds of dimension greater than $4$.
More information is needed.
Tolman and Weitsman used generating classes to compute the equivariant cohomology ring in case of a semifree action in [TW]. 
Their work gave us the idea for constructing necessary relations described in the present paper using information from generating classes.
Our proof was also motivated by the work of Goldin and Holm \cite{GH} where the Localization Theorem and dimensional reasoning were used.

\textbf{Organization}. In Section \ref{proof}, we prove our main result. Section \ref{examples} is devoted to several examples. 
Appendix \ref{generatorsfortoric} contains an algorithm for obtaining generating classes in the case of symplectic toric manifolds. This algorithm seems to be well known, however we could not find a good reference for it and therefore decided to include it in this paper for completeness.

\textbf{Acknowledgments.} The author is grateful to Tara Holm for suggesting this problem and for helpful conversations, 
and to the referees for their useful comments that allowed me to improve the exposition of the paper.
\section{Proof of Theorem \ref{main}}\label{proof}
Let a circle act on a manifold $M$ in a Hamiltonian fashion with isolated fixed points $p_1, \ldots, p_d$.
Let $\{a_p\}$ be a basis of $H_{S^1}^*(M; \bb{Q})$, satisfying condition ($\star$).
We want to show that if $f=(f_1, \ldots ,f_d) \in \oplus_{j=1}^d \,\bb{Q}[x]=H^*_{S^1}(M^{S^1})$ satisfies relations (\ref{relations}):
$$ \sum_{j=1}^d \frac{f_j\,a_p(p_j)}{e({p_j})} \in \bb{Q}[x],$$
for every fixed point $p$, then $f$
is in the image of injective, degree preserving map $\iota^*\colon H_{S^1}^*(M; \bb{Q}) \rightarrow H_{S^1}^*(M^{S^1}; \bb{Q})$. By abuse of notation we say such $f$ is an equivariant cohomology class of $M$. Recall that $e({p_j})$, the equivariant Euler class of the tangent bundle, evaluated at a fixed point $p_j$ is the product of the weights of the $S^1$ action on $T_{p_j}M$.
\begin{proof} Recall that $\bb{Q}[x]$ is a PID. Let $R$ be a submodule of $\oplus_{j=1}^d \,\bb{Q}[x]$ consisting of all $d$-tuples $f=(f_1, \ldots ,f_d)$ satisfying all of the above relations. As a submodule of a free module over PID, $R$ itself is free. Hamiltonian $S^1$-spaces are equivariantly formal, that is $H^*_{S^1}(M;\bb{Q}) \cong H^*(M;\bb{Q}) \otimes H^*(BS^1;\bb{Q})$ as modules. Therefore $\iota^*(H_{S^1}^*(M; \bb{Q}))$ is a free $\bb{Q}[x]$ submodule of $R \subset \oplus_{j=1}^d \,\bb{Q}[x]$. We already noticed that all the above relations are necessary. We show below that for any $k$ the number of generators of $(\,\iota^*(H_{S^1}^*(M; \bb{Q}))\,)_k$, the degree $k$ part of $\iota^*(H_{S^1}^*(M; \bb{Q}))$, is equal to the number of generators of $(R)_k$, the degree $k$ part of $R$. It then follows that $\iota^*(H_{S^1}^*(M; \bb{Q}))=R$ as needed.

We first analyze $\iota^*(H_{S^1}^*(M; \bb{Q}))$. The momentum map is a Morse function. Therefore the idex of a fixed point is well defined. Let $b_k$ be the number of fixed points of index $2k$. Then $d=\sum_{k=0}^n \, b_k$ is the  total number of fixed points.
By Theorem 1.3 of Frankel and Kirwan, we know that $b_k$ is also the $2k$-th Betti number of $M$. The fact $H^*_{S^1}(M;\bb{Q}) \cong H^*(M;\bb{Q}) \otimes H^*(BS^1;\bb{Q})$ implies that the equivariant Poincar\'{e} polynomial for $M$ is
$$ P^{S^1}_M(t)=P_M(t)P^{S^1}_{pt}(t)=(b_0 +b_1 t^2 + \ldots +b_n t^{2n})(1+ t^2+ t^4 +\ldots )= $$
$$=b_0 + (b_0 +b_1)t^2 + \ldots +(b_0 + b_1 + \ldots + b_k)t^{2k} + \ldots + d t^{2n} + d t^{2(n+1)} + \ldots .$$
Therefore $\iota^*(H_{S^1}^*(M; \bb{Q}))$ is a free $\bb{Q}[x]$ submodule of $R$, whose degree $k$ piece is a vector space over $\bb{Q}$ of dimension $(b_0 + b_1 + \ldots + b_k)$.

We now analyze and count the relations defining $R$.
For any  $f=(f_1, \ldots ,f_d) \in \oplus_{j=1}^d \,\bb{Q}[x]=H^*_{S^1}(M^{S^1})$ we denote by $K_j$ the degree of $f_j$ and by $r_{jk}\in n\bb{Q}$ its coefficients: $$f_j(x)=\sum_{k=0}^{K_j}r_{jk}x^k.$$
Then $r_{jk}$ are independent variables. Relations of type $$\sum_{j=0}^d\,s_j r_{jk}=0$$ 
for some constants $s_j$'s are called 
relations of degree $k$, as they involve the coefficients of $x^k$. Notice that if $f \in (\oplus_{j=1}^d \,\bb{Q}[x])_k$ is a homogeneous element of degree $k$ then it automatically satisfies all relations of degrees different then $k$. For any fixed point $p$ of index $2(k-1)$, a generating class $a_p$ associated with it assigns to each fixed point $p_j$ either $0$ or a homogeneous polynomial of degree $(k-1)$. 
Denote by $c^p_j$ the rational number satisfying $$\frac{a_p(p_j)}{e(p_j)}= c^{p}_{j}\,x^{k-1-n}.$$
If $f$ is an image of an equivariant cohomology class of $M$ then $f \cdot a_p$ is also. The Localization Theorem gives the relation 
$$ \int_M \, a_p f =\sum_{j=1}^d \frac{f_j\,a_p(p_j)}{e_{p_j}} \in \bb{Q}[x].$$
We may rewrite this in the following form:
\begin{align*}
 \int_M\, a_p f &=\sum_{j=1}^d \frac{f_j\,a_p(p_j)}{e(p_j)}\\
&=\sum_{j=1}^d f_j c^p_j\,x^{k-1-n}\\
&=\sum_{j=1}^d  c^p_j \left( \sum_{l=0}^{K_j}r_{jl}x^l \right)\,x^{k-1-n}\\
 &=\sum_{j=1}^d  c^{p}_{j} \left( \sum_{l=0}^{K_j}r_{jl}x^{k-1-n+l}\right) \in \bb{Q}[x].
 \end{align*}
Using the convention $r_{jl}=0$ for $l>K_j$, we can write
\begin{align*}
 \int_M\, a_p f &=\sum_{j=1}^d  c^{p}_{j} \left( \sum_{l=0}^{n-k}r_{jl}x^{k-1-n+l}\right)+\sum_{j=1}^d  c^{p}_{j} \left( \sum_{l=n-k+1}^{K_j}r_{jl}x^{k-1-n+l}\right)\\
                &= \sum_{l=0}^{n-k}\,\left( \sum_{j=1}^d  c^{p}_{j}  r_{jl}\right)x^{k-1-n+l}+\sum_{j=1}^d  c^{p}_{j} \left( \sum_{l=n-k+1}^{K_j}r_{jl}x^{k-1-n+l}\right).
 \end{align*}
The second component is an element of $\bb{Q}[x]$ as all the exponents of $x$ are nonnegative. 
Thus $\int_M\, a_p f$ is in $\bb{Q}[x]$ if and only if all the coefficients of $x$ in the first component 
(that is coefficients of negative powers of $x$) are $0$.
Therefore for any fixed point $p$ and any $l=0, \ldots, n-k$, where $2(k-1)$ is the index of $p$, we get the following linear relation of degree $l$:
$$ \sum_{j=1}^d  c^{p}_{j}r_{jl}=0.$$

Note that these relations are independent. 
We will show this by explicit computation. 
It is enough to show that for any $l$ all the relations of degree $l$ are independent, as relations of different degrees involve different subsets of variables $\{r_{jk}\}$.
Suppose that in some degree $l$ these relations in $r_{jl}$'s are not independent. That is, there are rational numbers $s_p$, not all zero, such that 
$$ \forall_{r_{jl}}\ \ \ \ 0=\sum_p s_p \left(\sum_{j=1}^d  c^{p}_{j}r_{jl} \right)=\sum_{j=1}^d \left(\sum_p s_pc^{p}_{j} \right)r_{jl}$$
As $r_{jl}$ are independent variables, we have
$ \sum_p s_pc^{p}_{j} =0$, for all $j=1, \ldots, d$.
Multiplying both sides by $e(p_j) x^{k-1-n}$ we obtain
$$ \sum_p s_p e(p_j) c^{p}_{j} x^{k-1-n}=0 .$$
Recall the definition of $c^{p}_{j}$ to notice that the above equation is equivalent to
$$\sum_p s_p a_p(p_j)=0.$$
That means $\sum_p s_p a_p$ vanishes on every fixed point and therefore is the $0$ class, 
although it is a nontrivial combination of classes $a_p$. 
This contradicts the independence of the generating classes $a_p$'s.

 Now we count the relations just constructed. As noted above, a fixed point of index $2(k-1)$
gives relations of degrees $0, \ldots, n-k$. Therefore a relation of degree $n-k$ is obtained from each fixed point of index $2(k-1)$ or less.
That means we get a relation of degree $k$ for each fixed point of index $2(n-k-1)$ or less, in total $$(b_0+ b_1 + \ldots +b_{n-k-1})$$
relations of degree $k$. The subspace of $(\oplus^d_{j=1} \bb{Q}[x])_k\cong \bb{Q}^d$ of elements satisfying all the relations of degree $k$ is of dimension $d-(b_0+ b_1 + \ldots +b_{n-k-1})$. Every homogeneous element $f \in (\oplus^d_{j=1} \bb{Q}[x])_k$ satisfying all the relations of degree $k$ also satisfies all the relations of other degrees (as coefficients of $x^l$ are $0$ for $l\neq k$). Moreover, the form of conditions (\ref{relations}) implies that for any $g \in \bb{Q}[x]$, $gf$ also satisfies all of relations (\ref{relations}).
Therefore the degree $k$ part of $R$ is the subspace of $(\oplus^d_{j=1}\bb{Q}[x])_k$ of elements satisfying all relations of degree $k$, and its dimension is $d-(b_0+ b_1 + \ldots +b_{n-k-1})$.
By the definition of $d$ and Poincar\'{e} duality,  
$$d-(b_0+ b_1 + \ldots +b_{n-k-1})=b_{n-k}+\ldots+b_n= b_0+ b_1 + \ldots + b_{k}.$$
This means that the degree $k$ part of $R$, $R_k$, is a vector space over $\bb{Q}$ of dimension $(b_0+ b_1 + \ldots + b_{k})$ containing a vector subspace $\iota^*(H_{S^1}^*(M; \bb{Q}))_k$, degree $k$ part of $\iota^*(H_{S^1}^*(M; \bb{Q}))$, of the same dimension. Therefore they must be equal. 
The two graded sumbodules: $\iota^*(H_{S^1}^*(M; \bb{Q}))$ and $R$, are equal in each degree. This implies 
$$\iota^*(H_{S^1}^*(M; \bb{Q}))=R.$$
\end{proof}

\section{Examples}\label{examples}
\begin{example}\label{toric1}
Consider the product of $\bb{C}P^2$ blown up at a point and $\bb{C}P^1$
$$\widetilde{\bb{C}P^2} \times \bb{C}P^1=\{([x_1:x_2],[y_0:y_1:y_2],[z_0:z_1])|\,x_1\,y_2-x_2\,y_1=0\},$$
and the following $T^3$ action on this space:
\footnotesize
$$(e^{iu},e^{iv},e^{iw})\cdot ([x_1:x_2][y_0:y_1:y_2][z_0:z_1])=([e^{iu}x_1:x_2][e^{iv}y_0:e^{iu}y_1:y_2][e^{iw}z_0:z_1]).$$
\normalsize
This is a symplectic toric manifold with 
moment map
$$ \mu([x_1:x_2][y_0:y_1:y_2][z_0:z_1])=
\left( \frac{|x_1|^2}{||x||^2}+\frac{|y_1|^2}{||y||^2},\frac{|y_0|^2}{||y||^2},
\frac{|z_0|^2}{||z||^2 } \right)$$
where $||x||^2=|x_1|^2+|x_2|^2$, and similarly for $||y||^2$ and $||z||^2$.
The moment polytope is shown in Figure \ref{fig:MomentPolyBlowUp}.
\begin{figure}[h]
\label{fig:MomentPolyBlowUp}
	\centering
		\includegraphics[width=0.60\textwidth]{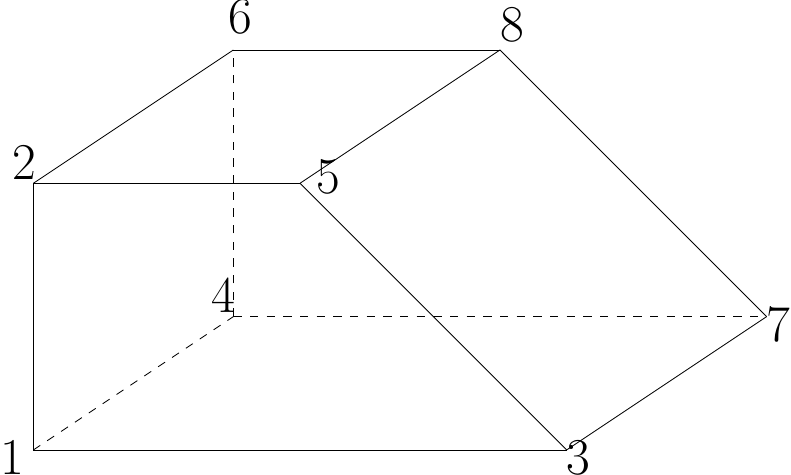}
	\caption{Moment polytope for $\widetilde{\bb{C}P^2} \times \bb{C}P^1$.}
\end{figure}
Using the algorithm from Appendix \ref{generatorsfortoric} we can compute generating classes for the equivariant cohomology with respect to $T$ action. 
They are presented in the table below.\\
\begin{tabular}[t]{c|ccccccccc}
class          & $v_1$ & $v_2$ &$v_3$ & $v_4$ & $v_5$    & $v_6$   & $v_7$  & $v_8$ & \\ \hline
$A_1$           & $1$   & $1$   & $1$  & $1$   & $1$      & $1$     & $1$    & $1$    & 
\includegraphics[width=0.21\textwidth]{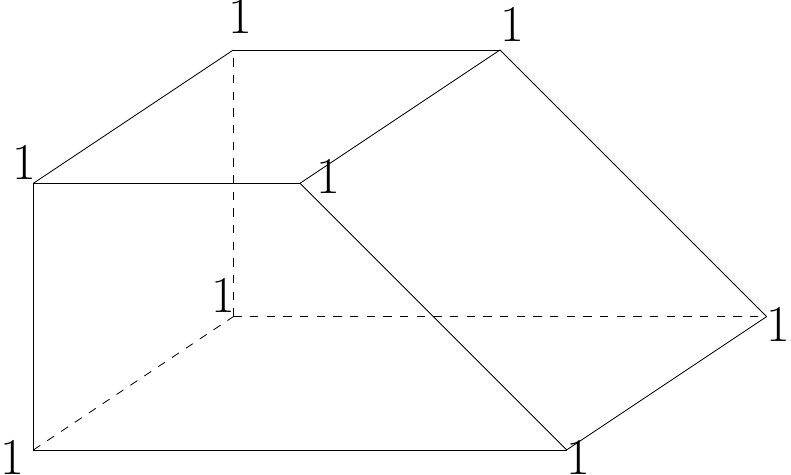} \\ 
$A_2$          & $0$   & $y$   & $0$  & $0$   & $y-x$    & $y$     & $0$    & $ y-x$ & 
\includegraphics[width=0.21\textwidth]{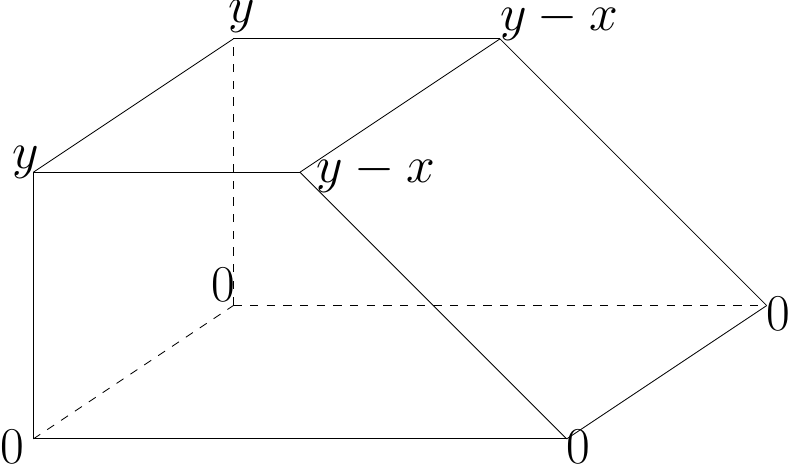} \\
$A_3$         & $0$   & $0$   & $x$  & $0$   & $x$      & $0$     & $x$    & $x$& 
\includegraphics[width=0.21\textwidth]{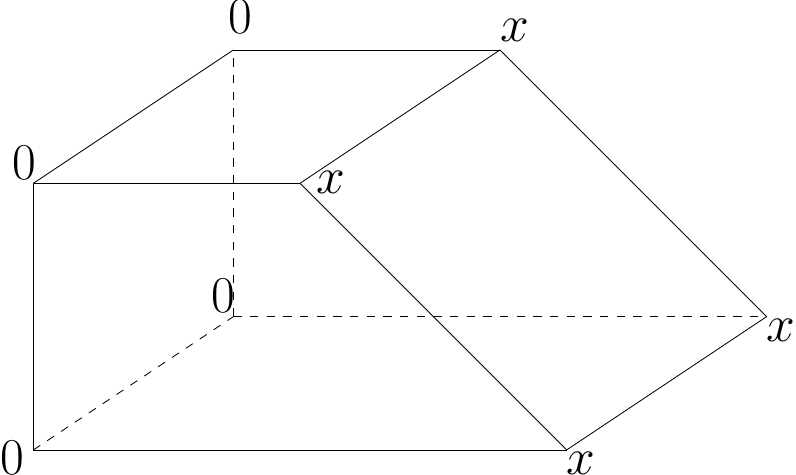}\\
$A_4$         & $0$   & $0$   & $0$  & $z$   & $0$      & $z$     & $z$    & $z$& 
\includegraphics[width=0.21\textwidth]{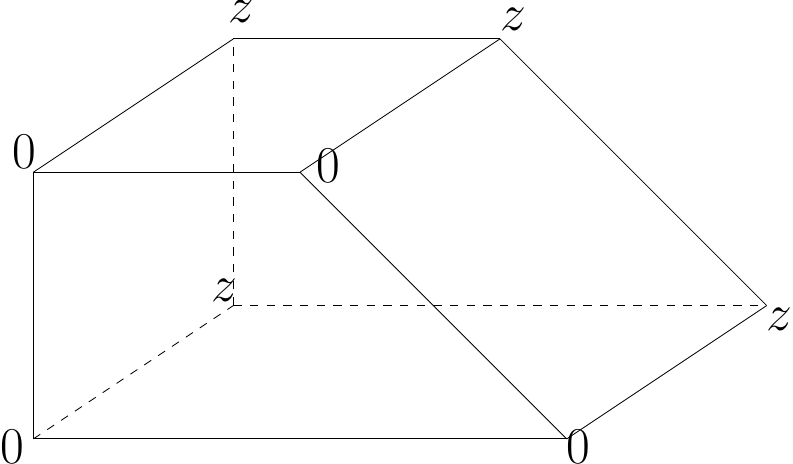}\\
$A_5$         & $0$   & $0$   & $0$  &  $0$  & $x(x-y)$ &$0$      & $0$    & $x(x-y)$& 
\includegraphics[width=0.21\textwidth]{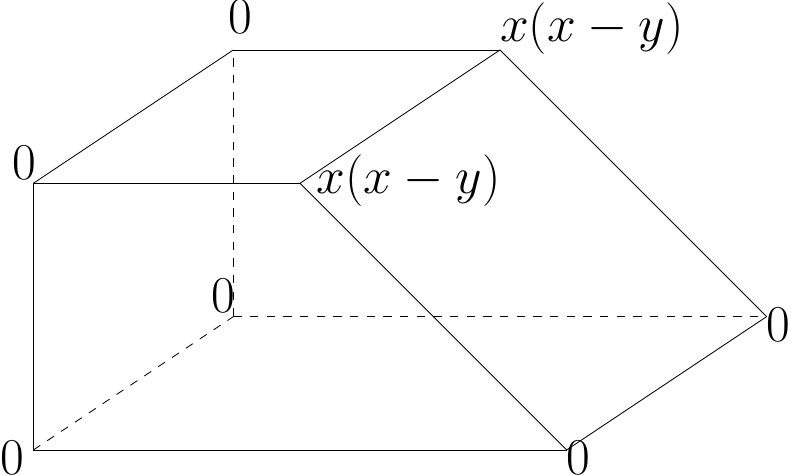}\\
$A_6$         & $0$   & $0$   & $0$  &  $0$  & $0$      & $yz$    & $0$    & $(y-x)z$& 
\includegraphics[width=0.21\textwidth]{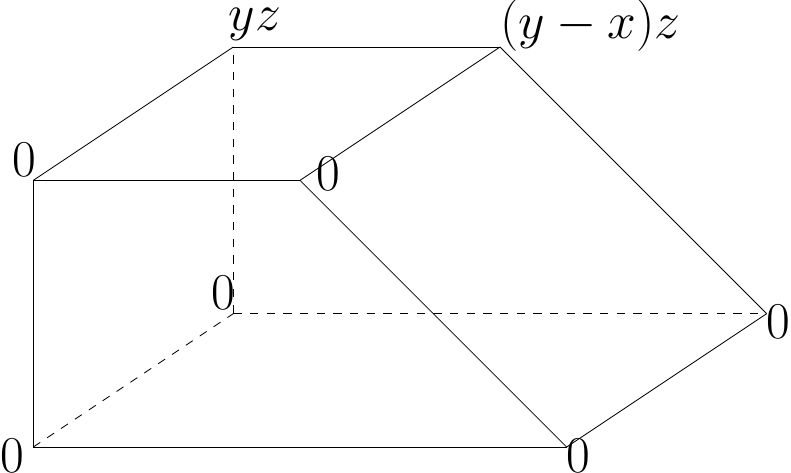}\\
$A_7$         & $0$   & $0$   & $0$  &  $0$  & $0$      & $0$     & $xz$   & $xz$& 
\includegraphics[width=0.21\textwidth]{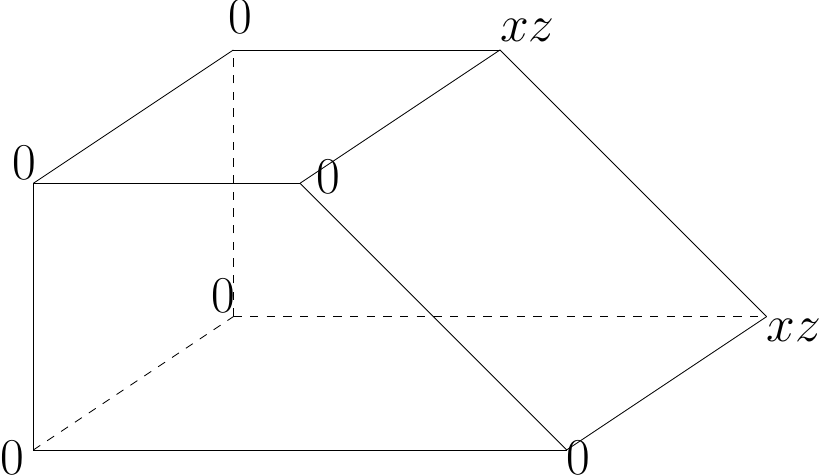}\\
$A_8$          & $0$   & $0$   & $0$  &   $0$ & $0$      &$0$      &$0$     & $xz(y-x)$& 
\includegraphics[width=0.21\textwidth]{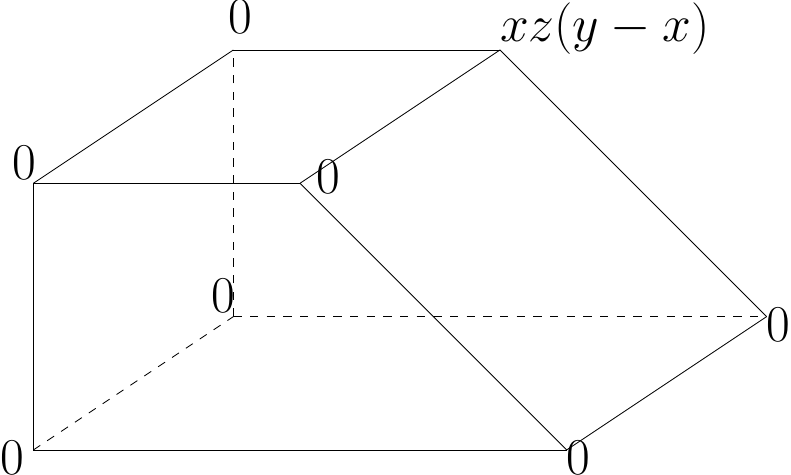}
\end{tabular}\\
\vskip 0.3cm
We want to compute equivariant cohomology with respect to the action of $S^1 \hookrightarrow T^3$ given by $u\rightarrow(u,2u,u)$. More precisely, our action is:
$$e^{iu} \cdot  ([x_1:x_2],[y_0:y_1:y_2],[z_0:z_1])=([e^{iu}x_1:x_2],[e^{i2u}y_0:e^{iu}y_1:y_2],[e^{iu}z_0:z_1]).$$
Note that we still have the same eight fixed points, namely:
 \vskip 0.2cm
\noindent
$v_1=([0:1],[0:0:1],[0:1]),$\\
$v_2=([0:1],[1:0:0],[0:1]),$\\
$v_3=([1:0],[0:1:0],[0:1]),$\\
$v_4=([0:1],[0:0:1],[1:0]),$\\
$v_5=([1:0],[1:0:0],[0:1]),$\\
$v_6=([0:1],[1:0:0],[1:0]),$\\
$v_7=([1:0],[0:1:0],[1:0]),$ and\\
$v_8=([1:0],[1:0:0],[1:0]).$\\
 \vskip 0.1cm
\noindent
The weights of this circle actions are:\\
\noindent
\begin{tabular}[t]{l|c|c}
fixed point & weights & index \\ \hline
$v_1$ & $u,2u,u$ & 0\\
$v_2$  & $u,-2u,u$ & 2 \\
$v_3$  & $-u,u,u$  & 2\\
$v_4$  & $u,2u,-u$  & 2\\
$v_5$  & $-u,-u,u$  & 4\\
$v_6$  & $u,-2u,-u$ & 4 \\
$v_7$  & $-u,u,-u$  & 4\\
$v_8$  & $-u,-u,-u$  & 6
\end{tabular}
\\
 \vskip 0.1cm
\noindent
We compute generating classes for the $S^1$ action from the classes for the $T$ action using the projection map $x\mapsto u$, $y \mapsto 2u$, $z \mapsto u$. 
They are presented in the table below, together with a row with $\frac{2u^3}{e(v_i)}$ that is useful for further computations.\\
\begin{tabular}[t]{c|cccccccc}
          & $v_1$ & $v_2$ &$v_3$ & $v_4$ & $v_5$    & $v_6$   & $v_7$  & $v_8$   \\ \hline
$\frac{2u^3}{e(v_i)}$ & 1&-1&-2&-1&2&1&2&-2\\ \hline        
$A_1$          & $1$   & $1$   & $1$  & $1$   & $1$      & $1$     & $1$    & $1$   \\ 
$A_2$          & $0$   & $2u$  & $0$  & $0$   & $u$      & $2u$    & $0$    & $ u$   
 \\
$A_3$          & $0$   & $0$   & $u$  & $0$   & $u$      & $0$     & $u$    & $u$ \\
$A_4$          & $0$   & $0$   & $0$  & $u$   & $0$      & $u$     & $u$    & $u$  \\
$A_5$          & $0$   & $0$   & $0$  &  $0$  & $u^2$    &$0$      & $0$    & $u^2$ \\
$A_6$          & $0$   & $0$   & $0$  &  $0$  & $0$      & $2u^2$  & $0$    & $u^2$ \\
$A_7$          & $0$   & $0$   & $0$  &  $0$  & $0$      & $0$     & $u^2$  & $u^2$ \\
$A_8$          & $0$   & $0$   & $0$  &   $0$ & $0$      &$0$      &$0$     & $u^3$ 
\end{tabular}
\\ 
 \noindent
 \vskip 0.1cm
We keep denoting by $f_j$ the restriction of $f$ to a fixed point $v_j$. The condition that 
$$  \sum_{j=1}^8\frac{f_j\,A_1}{e(V_j)}=\int_M f\,A_1 \in \bb{Q}[u],$$
implies that: \begin{displaymath}
         \begin{tabular}{lcl}
$\frac{f_1}{u^3}+\frac{-f_2}{u^3} +\frac{-2f_3}{u^3} +\frac{-f_4}{u^3} +\frac{ 2f_5}{u^3} +\frac{f_6}{u^3} +\frac{2f_7}{u^3} +\frac{-2f_8}{u^3}$  & $\in$ &$\bb{Q}[u]$. \\
\end{tabular}
\end{displaymath}
Thus\\
\begin{tabular}{lcll}
$f_1-f_2 -2f_3 -f_4 + 2f_5 +f_6 +2f_7 -2f_8 $ & $\in$ & $(u^3)$&$\bb{Q}[u]$, \\
\end{tabular}
\\ Similarly, using the class the $A_2$ we get\\
\begin{tabular}{rcll}
$-f_2 +f_5  +f_6 -f_8 $ & $\in$ & $(u^2)$&$\bb{Q}[u]$, \\
\end{tabular}
\\Other classes give:\\
\begin{tabular}{rcll}
$-f_3 +f_5  +f_7 -f_8 $   & $\in$ & $(u^2)$&$\bb{Q}[u]$, \\
$-f_4 +f_6  +2f_7 -2f_8 $ & $\in$ & $(u^2)$&$\bb{Q}[u]$, \\
$f_5 -f_8 $               & $\in$ & $(u)$&$\bb{Q}[u]$, \\
$2f_6 -2f_8 $             & $\in$ & $(u)$&$\bb{Q}[u]$, and \\
$f_7 -f_8 $               & $\in$ & $(u)$&$\bb{Q}[u]$.
\end{tabular}\\
Therefore $f=(f_1,\ldots, f_d)$ represents an equivariant cohomology class if and only if it satisfies:
\begin{itemize}
\item the degree $0$ relations:\\
$(f_i-f_j) \in (u)\bb{Q}[u]$, for every $i$ and $j$,
\item the degree $1$ relations:\\
\begin{tabular}{rcll}
$-f_3 +f_5  +f_7 -f_8 $   & $\in$ & $(u^2)$&$\bb{Q}[u]$ \\
$-f_2 +f_5  +f_6 -f_8 $   & $\in$ & $(u^2)$&$\bb{Q}[u]$ \\
$-f_4 +f_6  +2f_7 -2f_8 $ & $\in$ & $(u^2)$&$\bb{Q}[u]$ \\
$f_1 -f_2  -2f_3 +2f_5 $ & $\in$ & $(u^2)$&$\bb{Q}[u]$
\end{tabular}
\item the degree $2$ relation:\\
$f_1-f_2 -2f_3 -f_4 + 2f_5 +f_6 +2f_7 -2f_8 \in(u^3) \bb{Q}[u].$
\end{itemize}
\end{example}

\begin{example}\label{toric2}
 In the case of the specialization for a $T=T^n$ action on $M^{2n}$ (i.e. a symplectic toric manifold) to the action of some generic circle $S^1$ (i.e. with $M^{S^1}=M^T$), we can proceed using this simple algorithm. \\
\indent
The weights of $T$ action are easy to read from the moment polytope - they are just primitive integer vectors in the directions of the edges. 
To get the weights for our chosen $S^1$-action, we just need to use the appropriate projection $\pi: \lie{t}^* \rightarrow (\lie{s}^1)^*$. 
To find a basis of generating classes we first use the method from Appendix \ref{generatorsfortoric} with $\xi$ a generator of our $S^1$ to get a $T$-basis, and then use projection $\pi$. 
If the fixed points are $p_1, \ldots, p_d$, we denote by $a_1, \ldots, a_d$ the generating classes assigned to them and by $G_1, \ldots, G_d$ 
the faces of moment polytope that are the flow up faces of the corresponding fixed points. 
Recall that for any $v \in (\bb{Q}^n)^* \subset (\bb{R}^n)^*$ we denote by $prim(v)\in (\bb{Z}^n)^*$ the primitive integral vector in the direction of $v$.
Using this notation, and the construction from Appendix \ref{generatorsfortoric}, Theorem \ref{main} states that $f=(f_1, \ldots ,f_d) \in \oplus_{j=1}^d \,\bb{Q}[x]$ is an equivariant cohomology class of $M$ 
if and only if for any fixed point $p_l$ we have
$$ \sum_{j=1}^d \frac{f_j\,a_l(p_j)}{e(p_j)}= \sum_{\{j\,|\,p_j \in G_l\}} \frac{f_j\,\prod_{r}\pi(prim(r-p_j))}{e(p_j)}\in \bb{Q}[x],$$
where the product is taken over all vertices $r$ not in $G_l$ such that $r$ and $p_j$ 
are connected by an edge.
The equivariant Euler class $e(p_j)$ is the product of all weights at $p_j$, therefore, up to a multiplication by a rational constant, it is equal to
$$\prod_{r}\pi(prim(r-p_j)),$$
where the product is taken over all vertices $r$ connected to $p_j$. Thus the above condition is equivalent to
$$ \sum_{\{j\,|\,p_j \in G_l\}} \frac{f_j}{\prod_{r}\pi(prim(r-p_j))} \in \bb{Q}[x],$$
where the product is taken over all fixed points $r \in G_l$ that are connected with $p_j$ by an edge in $G_l$. 
\\
\indent
Consider, for example, vertex $v_3$ in the Example \ref{toric1} above. The face $G_3$ is the face spanned by $v_3,v_5,v_7,v_8$.
The weights at $v_3$ corresponding to edges that are in $G_3$ are $u,u$, for $v_5$: $u,-u$, for $v_7$: $-u,u$ and for $v_8$: $-u,-u$.
Therefore relation we get is:

$$ \frac{f_3}{u^2} +\frac{f_5}{-u^2} +\frac{f_7}{-u^2} +\frac{f_8}{u^2} \in \bb{Q}[u].$$

\noindent After clearing denominators, we obtain relation $f_3 -f_5 -f_7 +f_8 \in (u^2) \bb{Q}[u].$
\end{example}
\begin{example}\label{CSusage} Consider the following $T^2$ action on 
$M=\bb{C}P^4 $:
$$(e^{it},e^{is}) \cdot [z_0:z_1:z_2:z_3:z_4]=[z_0:e^{it}z_1:e^{i2t}z_2:e^{i3t}z_3:e^{i(t+s)}z_4].$$
This action has $5$ fixed points with the following weights:\\
\begin{tabular}[t]{c|c}
fixed point & weight  \\ \hline
$p_1=[1:0:0:0:0]$ & $x,2x,3x, x+y$ \\
$p_2=[0:1:0:0:0]$  & $-x,x,2x,y$  \\
$p_3=[0:0:1:0:0]$  & $-2x,-x,x, y-x$  \\
$p_4=[0:0:0:1:0]$  & $-3x,-2x,-x, y-2x$  \\
$p_5=[0:0:0:0:1]$  & $-x-y,-y,x-y,2x-y$  
\end{tabular}
\\ We want to find relations among $f_j$'s so that $f=(f_1, \ldots, f_5) \in \oplus_{j=1}^5 \bb{Q}[x,y]$ is in the image of
$$H^*_{T^2}(M)\hookrightarrow H^*_{T^2}(M^{T^2})= \oplus_{j=1}^5 \bb{Q}[x,y].$$
topological Schur Lemma, Theorem \ref{CSLemma}, this image is 
$$ \bigcap_H r^*_{M^H}(H^*_{T^2}(M^H)),$$
where intersection is taken over all codimension $1$ subtori $H$ which appear as isotropy groups of some elements of $M$ 
(that is $M^H \neq M^{T^2}$). We have chosen the identification $T^2\cong S^1 \times S^1$ with the first circle factor corresponding to $x$, 
and the second to $y$ variable in $H^*(BT)=\bb{Q}[x,y].$
In this example there are two relevant subgroups of $T^2$: $H_1=S^1 \times \{1\}\hookrightarrow T^2$ and 
$H_2=\{1\} \times S^1 \hookrightarrow T^2$.\\
In the first case, $M^{H_1}=\{[0:z_1:0:0:z_4]\}\cong \bb{C}P^1$ and $S^1\cong T^2/H_1$ acts on $M^{H_1}$ by 
$$e^{is}\cdot [0:z_1:0:0:z_4]=[z_0:z_1:0:0:e^{is}z_4].$$
There are
two fixed points: $p_2$ and $p_5$. We get the following relation in $\oplus_{j=1}^5 \bb{Q}[x,y]$ (see Example \ref{sphere}):
$$f_2-f_5 \in (y)\, \bb{Q}[x,y].$$
In the case of $H_2$ have that $M_2:= M^{H_2}=\{[z_0:z_1:z_2:z_3:0]\} \cong \bb{C}P^3$. Fixed points of this action are \\
\begin{tabular}[t]{l|c|c}
fixed point & weight  & index \\ \hline
$p_1=[1:0:0:0:0]$ & $x,2x,3x$ & $0$\\
$p_2=[0:1:0:0:0]$  & $-x,x,2x$ & $2$ \\
$p_3=[0:0:1:0:0]$  & $-2x,-x,x$  & $4$\\
$p_4=[0:0:0:1:0]$  & $-3x,-2x,-x$ & $6$   
\end{tabular}
\\The moment map is 
$$[z_0:z_1:z_2:z_3:0] \rightarrow -\frac{1}{2}(\frac{|z_1|^2}{\sum_{i=0}^3 |z_i|^2}+ \frac{2|z_2|^2}{\sum_{i=0}^3 |z_i|^2}+\frac{3|z_3|^2}{\sum_{i=0}^3 |z_i|^2}).$$
To find the relations we first need to compute generating classes. We easily get that:\\
\begin{tabular}[t]{c|cccc}
class          & $p_1$ & $p_2$ &$p_3$ & $p_4$ \\ \hline
$A_1$          & $1$   & $1$   & $1$  & $1$       \\ 
$A_2$          & $0$   & $-x$  & $ax$  & $bx$     \\
$A_3$          & $0$   & $0$   & $2x^2$  & $cx^2$   \\
$A_4$         & $0$   & $0$   & $0$  & $-6x^3$   
\end{tabular}\\
where $a,b,c$ are some parameters. Dimension reasons give $0=\int_{M_2} A_3= \frac{2x^2}{2x^3} + \frac{cx^2}{-6x^3}$ thus $c=6$. 
We would like to apply the Goldin-Tolman formula (Theorem 1.6 in \cite{GT}) to compute the values of other generating classes. 
Goldin and Tolman worked with a very special collection of generating classes, called the canonical classes.
The canonical class assigned to a fixed point $p$ needs to vanish at all other points of index less than or equal to the index $p$ 
(see comments below Theorem \ref{kirwanclasses}). 
In this particular example, all our fixed points are of different index and therefore the above classes are canonical classes in the sense of Goldin and Tolman.  
This allows us to apply Theorem 1.6 from [GT] and compute that $A_2 (p_4)=\frac{6a\,x}{4}$. 
Substituting this result into $0=\int_{M_2} A_2$ gives that $a=-2$ is the unique solution. Therefore generating classes are\\
\begin{tabular}[t]{c|cccc}
class          & $p_1$ & $p_2$ &$p_3$ & $p_4$ \\ \hline
$A_1$          & $1$   & $1$   & $1$  & $1$       \\ 
$A_2$          & $0$   & $-x$  & $-2x$  & $-3x$     \\
$A_3$          & $0$   & $0$   & $2x^2$  & $6x^2$   \\
$A_4$          & $0$   & $0$   & $0$  & $-6x^3$   
\end{tabular}
\\The relations we obtain in this way are\\
\begin{tabular}{lrcll}
$A_3:$  & $-f_3 -f_4 $   & $\in$ & $(x)$&$\bb{Q}[x,y]$, \\
$A_2:$  &$f_2-2f_3 +f_4 $   & $\in$ & $(x^2)$&$\bb{Q}[x,y]$, \\
$A_1:$  &$f_1-3f_2 +3f_3  -f_4  $ & $\in$ & $(x^3)$&$\bb{Q}[x,y]$.\\
\end{tabular}\\
Simplifying the relations and putting all the results together we get that $f=(f_1, \ldots, f_5) \in \oplus_{j=1}^5 \bb{Q}[x,y]$ is in the image of
$H^*_{T^2}(M)\hookrightarrow  \oplus_{j=1}^5 \bb{Q}[x,y]$ if and only if it satisfies\\
\begin{tabular}{rcll}
$f_i -f_j $               & $\in$ & $(x)$&$\bb{Q}[x,y]$, \\
$f_2-2f_3 +f_4 $             & $\in$ & $(x^2)$&$\bb{Q}[x,y]$, \\
$f_1-2f_2 +f_3 $             & $\in$ & $(x^2)$&$\bb{Q}[x,y]$, \\
$f_1-3f_2 +3f_3 -f_4 $               & $\in$ & $(x^3)$&$\bb{Q}[x,y].$
\end{tabular}\\
\end{example}
\begin{example} \label{SO5}
 Let $T\subset SO(5)$ be the maximal $2$-torus in $SO(5)$ and 
let $\mathcal{O}_{\lambda}$ be the coadjoint orbit of $SO(5)$ through a generic point $\lambda \in \lie{t}^*$, 
the dual of the Lie algebra $\lie{t}$ of $T$.
The torus $T$ acts on $\mathcal{O}_{\lambda}$ in a Hamiltonian fashion. 
We compute the $T/H$ equivariant cohomology of $M=\mathcal{O}_{\lambda} // H$, 
the symplectic reduction of $\mathcal{O}_{\lambda}$ by a circle $H \subset T$ 
fixing an $S^2$ in $\mathcal{O}_{\lambda}$ and
chosen so that the reduced space is a manifold.
The inclusion $\lie{h} \hookrightarrow \lie{t}$ induces the projection $\Pi_H :\lie{t}^* \rightarrow \lie{h}^*$. 
To obtain a moment map for action of $H$,  $\Phi_H : \mathcal{O}_{\lambda} \rightarrow \lie{h}^*$, we need to compose the moment map $\Phi_T$ 
for the $T$ action with this projection. We choose a regular value, $\mu$, of $\Phi_H$ and define 
$$M=\mathcal{O}_{\lambda} // H :=\Phi_H^{-1}(\mu) /H.$$
The residual action of $G:=T/H \cong S^1$ on $M$ is Hamiltonian and the moment map image can be identified with a slice of the moment polytope of $\mathcal{O}_{\lambda}$ 
presented in Figure \ref{fig:SO(5)Reduction}. 
\begin{figure}[htbp]
	\centering
		\includegraphics[width=1.0\textwidth]{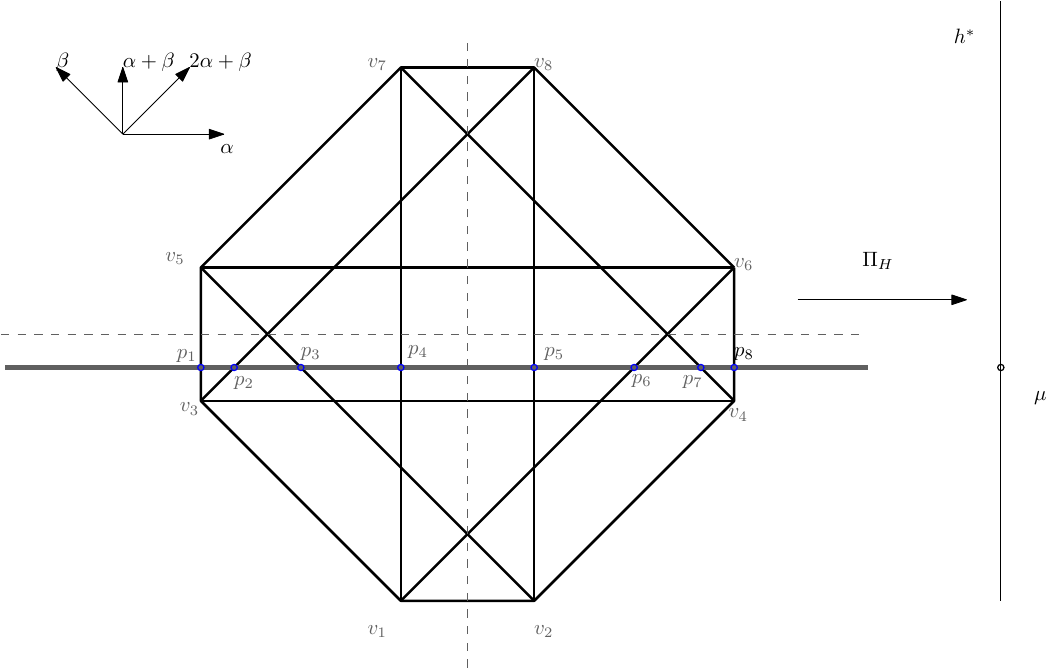}
		\caption{Moment Polytope for $T^2$ action on coadjoint orbit of  $SO(5)$ through a generic point. }
	\label{fig:SO(5)Reduction}
\end{figure}
We will compute the equivariant cohomology of $M$ with respect to this $T/H$ action.
This action has $8$ fixed points which we denote $p_1, \ldots, p_8$. For each $i$ there is a splitting of the torus $T=H\oplus H_i$ such that 
$H_i$ fixes a sphere $S_i^2\subset \mathcal{O}_{\lambda}$, and there is $q_i \in S_i^2$ such that $\Phi_T(q_i)=p_i$.
The residual $G=T/H$ action on $T_{p_i}M$ is isomorphic to the $H_i$ action on $N_{q_i} S_i^2$, the normal bundle to $S_i^2$
in $M$. To obtain the weights of $G$ action on $T_{p_i}M$, take the $T$ weights at the north or the south pole of  $S_i^2$
and 
compute their images under the projection $\lie{t}^* \rightarrow \lie{h}_i^*$. This projection is a map $\bb{Q}[\alpha, \beta] \rightarrow \bb{Q}[x]$ 
that sends $\alpha$ to $x$, and the weight assigned to $S_i$ to $0$. One of the $T$ weights will go to $0$ under this map. The three remaining weights
are the $G$ weights. Note that either pole will give the same result, 
as the $T$ weights differ by a multiple of the weight assigned to the edge representing $S_i$, and this weight vanishes on $\lie{h}_i$.
For our example we have\\
\begin{tabular}[t]{l|c|c}
fixed point & weight  & index \\ \hline
$p_1$ & $x,x,x$ & $0$\\
$p_2$  & $-x,x,2x$ & $2$ \\
$p_3$  & $-x,x,2x$  & $2$\\
$p_4$  & $-x,x,x$ & $2$  \\
$p_5$ & $-x,-x,x$ & $4$\\
$p_6$  & $-x,-2x,x$ & $4$ \\
$p_7$  & $-x,-2x,x$  & $4$\\
$p_8$  & $-x,-x,-x$ & $6$    
\end{tabular}\\
\vskip0.01cm
\noindent
We will use generating classes for the $T$ action on the whole coadjoint orbit to obtain generating classes for the $G$-equivariant cohomology of $M$.
Note that there will be $8$ generating classes for $H^*_G(M)$.
Let $A_j$ denote a canonical class for the $T$ action associated to a fixed point $v_j$, for $j=1,...,8$, with $\alpha$, $\beta$ as in Figure \ref{fig:SO(5)Reduction}.
Figure \ref{SO5Chart} presents their values in a chart, while 
Figure \ref{CanClassesSO(5)} presents them graphically.\\
\begin{figure}[h]
 \includegraphics[width=1.0\textwidth]{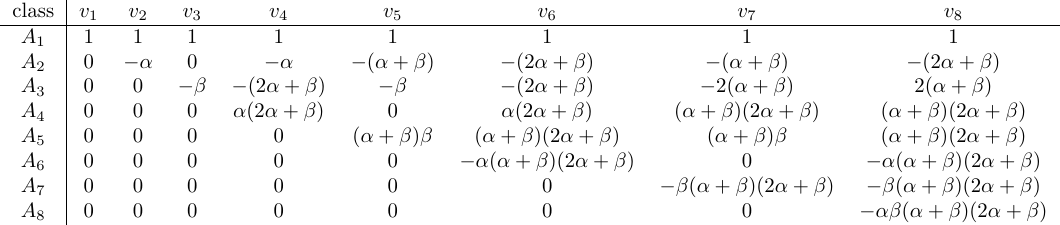}
\caption{The values of generating classes for $T^2$ action on coadjoint orbit of $SO(5)$.}
\label{SO5Chart}
\end{figure} 
\begin{figure}
	\includegraphics[width=1.0\textwidth]{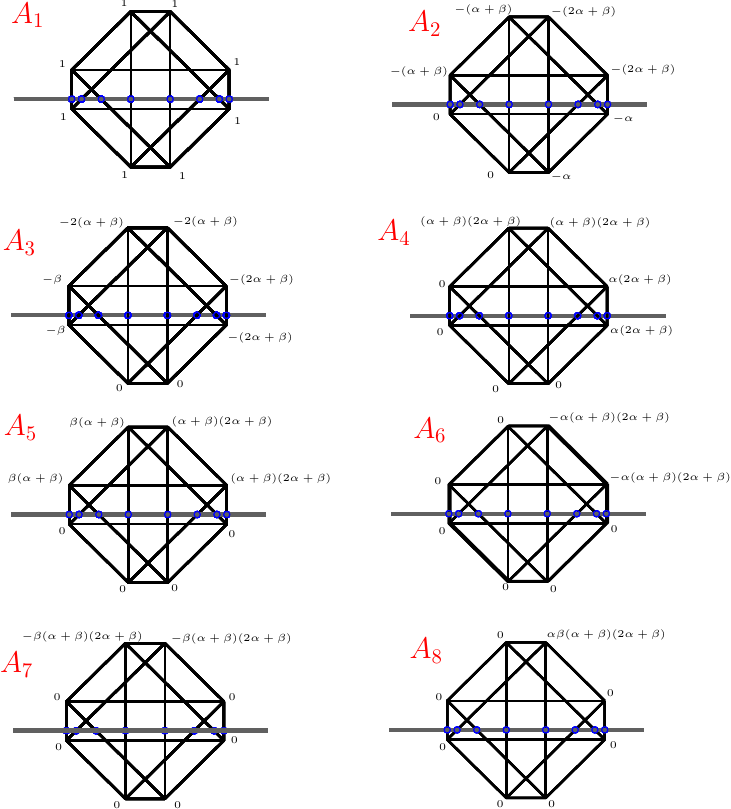}
	\caption{The generating classes for $T^2$ action on coadjoint orbit of $SO(5)$.}
	\label{CanClassesSO(5)}
\end{figure}
There is a surjective map $\kappa: H^*_{T}(\mathcal{O}_{\lambda}) \rightarrow H^*_{G}(M)$, and therefore 
generating classes for $G$ action on $M$ are images of some $\bb{Q}[\alpha, \beta]$-linear combinations of $A_j$'s. 
To compute the value of $\kappa(A)$ at $p_i$ for $A \in H^*_{T}(\mathcal{O}_{\lambda})$, 
take the value of $A$ at north or south pole of $S_i$ and 
compute its image under the map $\bb{Q}[\alpha, \beta] \rightarrow \bb{Q}[x]$ 
that sends $\alpha$ to $x$, and the weight assigned to $S_i$ to $0$.
 The kernel of map $\kappa$ was described by Tolman and Weitsman 
in \cite{TW:abelianquotient}. Their result implies that for $i=5, \ldots 8$, $\kappa(A_i)=0$. 
Therefore the image of $\kappa$ is generated over 
$\bb{Q}[x]$ by $\kappa(A_1), \ldots,\kappa(A_4),\kappa(\beta A_1), \ldots, \kappa(\beta A_4)$.
Moreover, these $8$ classes are $\bb{Q}[x]$-linearly independent.
They are not Kirwan classes as described in Theorem \ref{kirwanclasses}. 
However they satisfy all the requirements needed to apply our Theorem, that is condition ($\star$). 
They are in bijection with the fixed points and 
a class corresponding to a fixed point 
of index $2k$ evaluated at any fixed point is $0$ or a homogenous polynomial of degree $k$.
We present them, together with the Euler class, in the following table:\\
\begin{center}
\begin{tabular}[t]{l|ccccccccc}
class           & $p_1$ & $p_2$ &$p_3$ & $p_4$ & $p_5$    & $p_6$   & $p_7$  & $p_8$\\ \hline
$a_1=\kappa(A_1)$           & $1$   & $1$   & $1$  & $1$   & $1$      & $1$     & $1$    & $1$     \\ 
$a_2=\kappa(A_2)$           & $0$   & $0$  & $-x$  & $0$   & $-x$      & $0$    & $-x$    & $ -x$  \\
$a_3=\kappa(A_3)$           & $x$   & $2x$   & $0$  & $0$   & $x$      & $0$     & $-2x$    & $-x$\\
$a_4=\kappa(\beta A_1)$           & $-x$   & $-2x$   & $0$  & $-x$   & $-x$      & $-2x$     & $0$    & $-x$\\
$a_5=\kappa(A_4)$           & $0$   & $0$   & $0$  &  $0$  & $0$      & $0$     & $2x^2$  & $x^2$\\
$a_6=\kappa(\beta A_2)$           & $0$   & $0$   & $0$  &  $0$  & $x^2$    &$0$      & $0$    & $x^2$\\
$a_7=\kappa(\beta A_3)$           & $-x^2$   & $-4x^2$   & $0$  &  $0$  & $0$      & $0$  & $0$    & $x^2$\\
$a_8=\kappa(\beta A_4)$           & $0$   & $0$   & $0$  &   $0$ & $0$      &$0$      &$0$     & $-x^3$\\
Euler class    &$x^3$&  $-2x^3$&  $-2x^3$& $-x^3$&  $x^3$&  $2x^3$&  $2x^3$&  $-x^3$
\end{tabular}\\
\end{center}
\vskip 0.1cm
Therefore we get the following relations on $f=(f_1, \ldots, f_8)$:\\
\begin{tabular}{lrcll}
$a_1:$  & $2f_1-f_2-f_3 -2f_4 +2f_5+f_6+f_7-2f_8 $   & $\in$ & $(x^3)$&$\bb{Q}[x]$ \\
$a_2:$  &$f_3 -2f_5-f_7+2f_8 $   & $\in$ & $(x^2)$&$\bb{Q}[x]$ \\
$a_3:$  &$2f_1-2f_2 -2f_7+2f_8 $ & $\in$ & $(x^2)$&$\bb{Q}[x]$ \\
$a_4:$  &$-2f_1+2f_2 +2f_4 -2f_5-2f_6+2f_8   $ & $\in$ & $(x^2)$&$\bb{Q}[x]$ \\
$a_5:$  &$2f_7 -2f_8 $               & $\in$ & $(x)$&$\bb{Q}[x]$ \\
$a_6:$  &$2f_5 -2f_8$             & $\in$ & $(x)$&$\bb{Q}[x]$\\ 
$a_7:$  &$-2f_1+4f_2-2f_8 $               & $\in$ & $(x)$&$\bb{Q}[x]$ 
\end{tabular}\\
\vskip 0.1cm
Simplifying these relations we can put them in the following form:
\begin{itemize}
\item the degree $0$ relations:\\
All $(f_i-f_j) \in (x)\bb{Q}[x]$,
\item the degree $1$ relations:\\
\begin{tabular}{rcll}
$f_3 -2f_5-f_7+2f_8 $   & $\in$ & $(x^2)$&$\bb{Q}[x]$ \\
$f_1-f_2 -f_7+f_8 $   & $\in$ & $(x^2)$&$\bb{Q}[x]$ \\
$-f_3+3f_4-f_5-f_6 $ & $\in$ & $(x^2)$&$\bb{Q}[x]$ \\
$f_2  -f_3 -f_6 +f_7 $ & $\in$ & $(x^2)$&$\bb{Q}[x]$
\end{tabular}
\item the degree $2$ relation:\\
$2f_1-f_2 -f_3 -2f_4 + 2f_5 +f_6 +f_7 -2f_8 \in(x^3) \bb{Q}[x].$
\end{itemize}
\end{example}


\appendix
\section{Generating classes for Symplectic Toric Manifolds and their specializations.}\label{generatorsfortoric}
\begin{figure}[h]
\centering
	\includegraphics[width=0.70\textwidth]{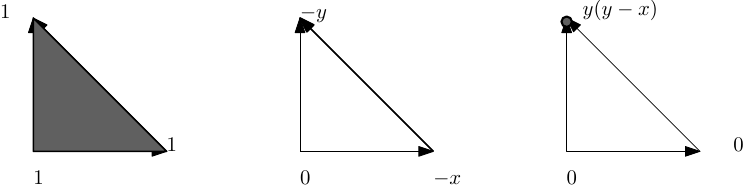} \\
	\caption{Generating classes for $\bb{C}P^2$.}
	\label{fig:CP2}
\end{figure}

\begin{figure}[h]
	\centering
		\includegraphics[width=0.62\textwidth]{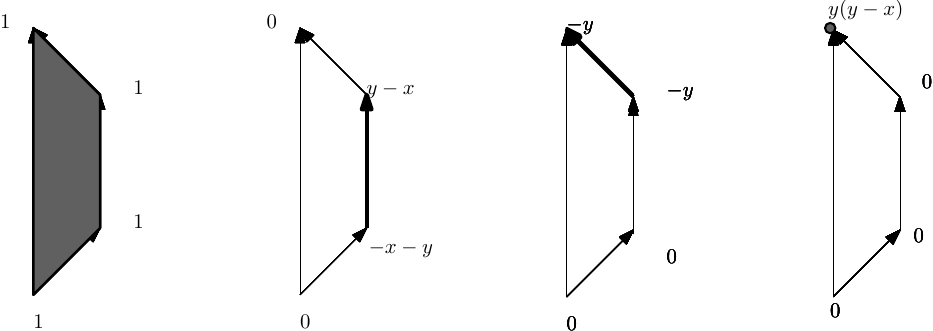}
			\\
			\caption{The basis of the equivariant cohomology ring given by the above algorithm.}
			\label{fig:Not unique generating classes}
	\end{figure}

\begin{figure}[h]
	
	\centering
		\includegraphics[width=0.62\textwidth]{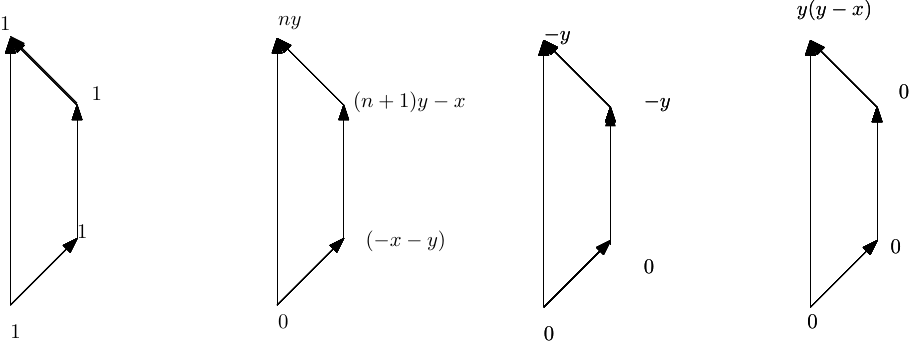}
\\
	\caption{Different basis of the equivariant cohomology ring.}
	\label{fig:differentbasis}
\end{figure}

A {\bf symplectic toric manifold} is a connected symplectic manifold $(M,\omega)$ equipped with an effective Hamiltonian action of a torus $T$ of dimension $\dim T=\frac 1 2 \dim_{\bb{R}} M$.
Let $M^{2n}$ be a compact symplectic toric manifold with a momentum map image a Delzant polytope $\Phi(M)=P \subset \lie{t}^*$. In particular the polytope $P$ is simple, rational and smooth.
The Lie algebra dual, $\lie{t}^*$, is isomorphic to $\bb{R}^n$, though not canonically. 
One of the conventions is to identify $S^1$ with $\bb{R}/\bb{Z}$. Then the exponential map $Lie(S^1)\cong \bb{R} \rightarrow S^1$ is of the form $t \rightarrow e^{2\pi i t}$.
With this identification, the function 
$$  \bb{C} \ni z \rightarrow -\pi \,k\,|z|^2 \in \bb{R}\cong Lie (S^1)$$
is a momentum map for the $S^1$ action on $(\bb{C}, \omega_{standard})$ by rotation with weight $k$.
Identifying $\lie{t}^*$ with $\bb{R}^n$ using above convention allows us to think of $P$ as a Delzant polytope in $\bb{R}^n$.
Denote by $M_1$ the union of all $T$-orbits of dimension $1$.
The closures of the connected components of $M_1$ are spheres, called the isotropy spheres.
Denote by $V$ the vertices of $P$, and by $E$ the $1$-dimensional faces of $P$, also called edges.
Vertices correspond to the fixed points of the torus action, while edges correspond to the isotropy spheres.
Fix a generic $\xi \in \bb{R}^n$, so that for any $p,q \in V$ we have $\langle p,\xi\rangle \, \neq\,\langle q,\xi\rangle$.
Orient the edges so that 
$\langle i(e), \xi \rangle \,\,< \,\, \langle t(e), \xi \rangle $
for any edge~$e$, where $i(e),\,t(e)$ are the initial and the terminal points of $e$. 
Let $w_{i(e)}(e)=-w_{t(e)}(e)$ denote the isotropy weights of the $T$ action on the tangent spaces to the isotropy sphere $\Phi \inv(e)$, $T_{\Phi \inv (i(e))}{\Phi \inv (e)}$ and $T_{\Phi \inv (t(e))}(\Phi \inv (e))$ respectively.
Note that $w_{i(e)}(e)$ is the primitive integral vector in the direction of $\vec{e}$. We denote it by $prim(\vec{e})$.
For any $p \in V$ let $G_p$ denote the smallest face containing $p$ and all points $q \in V$ 
with $\langle p, \xi\rangle\,<\,\langle q, \xi\rangle$ which are connected with $p$ by an edge . We will call $G_p$ the {\bf flow up face} for $p$. 
We define the class $a_p \in H_{S^1}^*(M^{S^1})$ by
\begin{displaymath} a_p(q)=
\begin{cases} 0\ & for\ q \in V\setminus G_p\\
\prod_{r}\,prim(r-q)\ & for\ q \in  G_p
\end{cases}
\end{displaymath}
where the product is taken over all $r \in V\setminus G_p$ such that $r$ and $q$ 
are connected by an edge of $P$. 
We use the convention that the empty product is $1$. 
If $k$ edges terminate at $p$ then the $n-k$ edges starting from $p$ belong to the face $G_p$ 
(as the polytope is simple, exactly $n$ edges meet at each vertex).
The smoothness of $P$ implies that these $n-k$ edges span an $(n-k)$ affine hyperplane $H_p$ of $\bb{R}^n$ and the face $G_p$ is the intersection $G_p=P\cap H_p$.
Moreover, it also implies that for any $q \in G_p$ there are $n-k$ edges meeting at $q$ that are contained in the face $G_p$
and $k$ edges connecting $q$ to vertices outside the face $G_p$.
Therefore the class $a_p$ assigns to each fixed point $0$ or a homogeneous polynomial of degree $k$. 
Such classes satisfy the GKM conditions and thus 
are in the image of the equivariant cohomolgy of $M$. 
The class $a_p$ constructed this way is the canonical equivariant extension (see \cite{LS}, Corollary 3.5) 
of the cohomology class Poincar\'{e} dual to the submanifold of $M$ mapping to the face $G_p$.
These two facts can be proved using the notion of the axial function introduced in \cite{GZ}. 
The classes we have just defined are also linearly independent, 
which follows easily from the fact that $a_p$ can be nonzero only at vertices $q$ greater or equal to $p$ in the partial order given by the orientation of edges. For example the classes presented in Figure \ref{fig:CP2} form a basis of generating classes for $\bb{C}P^2$ .
Recall that a basis of generating classes does not need to be unique. For example Figures \ref{fig:Not unique generating classes} and \ref{fig:differentbasis} present two different bases of generating classes for the equivariant cohomology ring of Hirzebruch surface $F_2$. The one on Figure \ref{fig:Not unique generating classes} is obtained using above algorithm.
	\\

This algorithm is also very useful while dealing with specialization, that is while restricting a toric $T$ action on $M$ to an action of some subtorus $S \hookrightarrow T$ (not necessarily a circle). As explained in the introduction, if $S$ is generic then $M^T=M^S$. To find a basis of generating classes for $H_S^*(M)$ we don't need to know the polytope for the $T$-action, nor the weights of the $T$ action. It is enough to know the isotropy weights of the $S$ action, 
the fact that this action is a specialization of some toric action and positions of the isotropy spheres for that toric action.
These weights are just projections of $T$ weights under $pr: \lie{t}^* \rightarrow \lie{s}^* $. 
That is the $S$ weight on edge $e$ is $pr \,(prim(t(e)-i(e))\,)$.
The positions of isotropy spheres for the toric action allow us to find the flow up face $G_p$  for any fixed point $p$.
The above algorithm gives that
\begin{displaymath} a_p(q)= 
\begin{cases} 0\ & for\ q \in V\setminus G_p\\
 \prod_{r}\,pr\,(prim(r-q)\,)\ & for\ q \in  G_p
\end{cases}
\end{displaymath}
where the product is taken over all $r \in V\setminus G_p$ such that $r$ and $q$ 
are connected by an isotropy sphere. If $S$ is a circle we may apply Theorem \ref{main} to obtain all relations needed to describe $\iota^*(H^*_{S^1}(M))$.
This gives us a method for computing equivariant cohomology for these circle action 
that could be extended a toric action. If $S$ is of bigger dimension, we need to apply Theorem \ref{CSLemma} together with Theorem \ref{main}.

\end{document}